\documentclass[11pt]{article}
\usepackage{graphicx, colordvi}
\usepackage{dsfont}
\usepackage{epsfig}
\usepackage{mathrsfs}
\usepackage{amssymb,amsfonts,amsmath,amsthm,cite,color}
\usepackage{dsfont}
\usepackage{epsfig}
\usepackage{mathrsfs}
\usepackage{longtable}
\usepackage{enumerate}
\usepackage{geometry}
\usepackage{hyperref}
\hypersetup{
colorlinks=true,
linkcolor=cyan,
filecolor=blue,
urlcolor=red,
citecolor=green,
}

\newtheorem{theo}{Theorem}

\newtheorem{coro}[theo]{Corollary}

\newtheorem{lem}[theo]{Lemma}

\makeatletter \@addtoreset{equation}{section}
\@addtoreset{theo}{section} \makeatother

\newcommand{\bN} { {\mathbb{N}}}

\newcommand{\bZ} { {\mathbb{Z}}}

\newcommand{\bK} { {\mathbb{K}}}

\DeclareMathOperator{\ann}{ann}

\newcommand{\la} { {\langle}}
\newcommand{\ra} { {\rangle}}
\newcommand{\re}{\noindent{\bfseries Remark. }}

\begin{document}
\begin{center}
 {\large \bf Power-Partible Reduction and Congruences for Ap\'ery Numbers}
\end{center}

\begin{center}
{  Rong-Hua Wang}$^{1}$ and {Michael X.X. Zhong}$^{2}$

   $^1$School of Mathematical Sciences\\
   Tiangong University \\
   Tianjin 300387, P.R. China\\
   wangronghua@tiangong.edu.cn \\[10pt]

   $^2$School of Science\\
   Tianjin University of Technology \\
   Tianjin 300384, P.R. China\\
   zhong.m@tjut.edu.cn
\end{center}

\begin{abstract}
In this paper, we introduce the power-partible reduction for holonomic (or, P-recursive) sequences and apply it to obtain a series of congruences for Ap\'ery numbers $A_k$.
In particular, we prove that, for any $r\in\bN$,
there exists an integer $\tilde{c}_r$ such that
\begin{equation*}
\sum_{k=0}^{p-1}(2k+1)^{2r+1}A_k\equiv \tilde{c}_r p \pmod {p^3}
\end{equation*}
holds for any prime $p>3$.
\end{abstract}



\section{Introduction}

In the 1990s, Wilf and Zeilberger \cite{WZ1990,Zeilberger1990c, Zeilberger1991} developed the mechanic proof theory
of combinatorial identities (i.e., the WZ theory).
Zeilberger's algorithm, or the method of creative telescoping, is the core of the WZ theory.
The reduction-based algorithms are a new class of algorithms for creative telescoping, which separate the calculation of telescopers from costly computation of certificates \cite{Chen2019}.
In particular, a polynomial reduction process for hypergeometric terms was introduced in 2015 by Chen et al. \cite{CHKL2015} in the development of the modified Abramov--Petkov{\v{s}}ek reduction.

In 2021, Hou, Mu and Zeilberger \cite{HouMuZeil2021} presented another polynomial reduction process, which avoids the multiplicative decomposition needed in  \cite{CHKL2015}.
This polynomial reduction was employed by Hou and Li \cite{HouLi2021} to derive new hypergeometric identities.
A $q$-analogue of the Hou--Mu--Zeilberger reduction has also been introduced to prove and discover $q$-identities automatically \cite{WZ2022a}.
Recently, the authors generalized the Hou--Mu--Zeilberger reduction to the holonomic case in \cite{WZ2022b}.
This provides an algorithmic way to prove and discover new multi-sum identities.

Given a polynomial $Q(k)$ and a hypergeometric term $t_k$, the Hou--Mu--Zeilberger reduction \cite{HouMuZeil2021} rewrites $Q(k)t_k$ as
\begin{equation*}Q(k)t_k=\Delta (T_k)+\tilde{q}(k)t_k,\end{equation*}
where $\Delta$ is the difference operator, $T_k$ is a hypergeometric term and $\tilde{q}(k)$ is a reduced polynomial such that $\deg \tilde{q}(k)\leq \deg Q(k)$.
A key observation in \cite{HouMuZeil2021} is that when the hypergeometric term $t_k$ satisfies certain symmetry conditions about a constant $\gamma$, the reduced polynomial $\tilde{q}(k)$ contains only odd (or, even) powers of $(k-\gamma)$ if the same goes for $Q(k)$.
We call such a hypergeometric term $t_k$ \emph{power-partible} with respect to $\gamma$ under the reduction.

In this paper, we investigate the power-partibility of holonomic sequences under the polynomial reduction developed in \cite{WZ2022b}.
As applications, we deduce new congruences involving Ap\'ery numbers $A_k$.

Ap\'ery numbers, arising from Ap\'ery's \cite{Poorten1979} proof of the irrationality of $\zeta{(3)}$, are defined as
\begin{equation*}
A_n=\sum_{k=0}^{n}\binom{n}{k}^2\binom{n+k}{k}^2.
\end{equation*}
In \cite{Sun2012}, Sun showed that for any prime $p>3$,
\begin{equation*}
\sum_{k=0}^{p-1}(2k+1)A_k\equiv p \pmod {p^3},
\end{equation*}
and conjectured that
\begin{equation}\label{eq:(2k+1)Sun-Guo-Zeng}
\sum_{k=0}^{p-1}(2k+1)(-1)^kA_k\equiv p\left(\dfrac{p}{3}\right) \pmod {p^3},
\end{equation}
which was later proved by Guo and Zeng \cite{GuoZeng2012}.
Here $\left(\dfrac{\bullet}{p}\right)$ denotes the Legendre symbol.

In \cite{Sun2016}, Sun further conjectured that for each $r\in\bN$ and prime $p>3$, there is a $p$-adic integer $c_r$ only depending on $r$ such that
\begin{equation}\label{eq: Sun-Conj-Alternating}
\sum_{k=0}^{p-1}(2k+1)^{2r+1}(-1)^kA_k\equiv c_r p \left(\dfrac{p}{3}\right) \pmod {p^3},
\end{equation}
which was confirmed very recently by Xia and Sun \cite{XiaSun2022}.
For more interesting arithmetic properties of $A_k$, see \cite{Gessel1982,GuoZeng2012,Pan2014,Sun2012,Sun2016,XiaSun2022}.

Using the power-partible reduction for holonomic sequences developed in Section 2, we prove that for each $r\in\bN$,
there exists an integer $\tilde{c}_r$ such that
\begin{equation*}
\sum_{k=0}^{p-1}(2k+1)^{2r+1}A_k\equiv \tilde{c}_r p \pmod {p^3}
\end{equation*}
holds for any prime $p>3$.
The proof also provides a constructive algorithm to calculate $\tilde{c}_r$ for any given positive integer $r$.

The rest of the paper is organized as follows.
A criterion (Theorem \ref{th:power reduce v2}) on the power-partibility of holonomic sequences is presented in Section 2, which is a generalization and unification of Hou--Mu--Zeilberger's analogous results for hypergeometric terms.
New series of congruences for Ap\'ery numbers are obtained in Section 3.

\section{Polynomial reduction and power-partibility}

Let $\bK$ be a field of characteristic $0$ and $\bK[k]$ the polynomial ring over $\bK$.
The set of \emph{annihilators} of a sequence $F(k)$ is denoted by
\begin{equation}\label{eq:Annihilator}
\ann F(k):=\left\{L=\sum_{i=0}^{J}a_i(k)\sigma^i\in \bK[k][\sigma]\mid L(F(k))=0\right\},
\end{equation}
where $J\in\bN=\{0,1,2,\ldots\}$ and $\sigma$ is the shift operator (that is, $\sigma F(k)=F(k+1)$).
A sequence $F(k)$ is said to be \emph{holonomic} (or, \emph{P-recursive}) if $\ann F(k)\neq \{0\}$.
We call $J$ in Eq.~\eqref{eq:Annihilator} the \emph{order} of $L$ if $a_J(k)\neq 0$,
and the minimum order of all $L\in\ann F(k)\setminus\{0\}$ is called the \emph{order} of $F(k)$.

A holonomic sequence $F(k)$ of order $J$ is called \emph{summable} if one can write
\begin{equation*}
F(k)=\Delta\left(\sum_{i=0}^{J-1}u_i(k)F(k+i)\right)
\end{equation*}
for some rational functions $u_i(k)\in\bK(k)$, where $\Delta=\sigma-1$ denotes the difference operator.
Given a holonomic sequence $F(k)$, one can always construct polynomials $q(k)$ such that $q(k)F(k)$ is summable.

For any operator $L=\sum_{i=0}^{J}a_i(k)\sigma^i$ with $a_i(k)\in\bK[k]$,
the \emph{adjoint} of $L$ is defined by
\begin{equation*}\label{eq:L}
L^{\ast}=\sum_{i=0}^{J}\sigma^{-i}a_i(k).
\end{equation*}
Namely, for any polynomial $x(k)\in\bK[k]$, we have
\begin{equation*}\label{eq:B}
L^{\ast}(x(k))=\sum_{i=0}^{J}a_i(k-i)x(k-i).
\end{equation*}

\begin{theo}[van der Hoeven \cite{Hoeven2018}]
Suppose that $(F(k))_{k=0}^{\infty}$ is a holonomic sequence and that $L=\sum_{i=0}^{J}a_i(k)\sigma^i$ $\in\ann F(k)$ with $a_J(k)\neq 0$.
Then for any $x(k)\in\bK[k]$,
\begin{equation}\label{eq: summable}
L^{\ast}(x(k))F(k)=\Delta\left(-\sum_{i=0}^{J-1}u_i(k)F(k+i)\right),
\end{equation}
where
\begin{equation}\label{eq:u}
u_i(k)=\sum_{j=1}^{J-i}a_{i+j}(k-j)x(k-j),\quad i=0,1,2,\ldots,J-1.
\end{equation}
\end{theo}

\re
Summing over $k$ from $0$ to $n-1$ on both sides of Eq.~\eqref{eq: summable}, we obtain
\begin{equation}\label{eq:rec ann}
\sum_{k=0}^{n-1}L^{\ast}(x(k)) F(k)
=\left(\sum_{i=0}^{J-1}u_i(0)F(i)\right)
-\left(\sum_{i=0}^{J-1}u_i(n)F(n+i)\right).
\end{equation}

The set
\begin{equation*}\label{eq:SL}
S_{L}=\{L^{\ast}(x(k))\mid x(k)\in\bK[k]\}
\end{equation*}
is called the \emph{difference space} corresponding to $L$.
Let $[p(k)]_L=p(k)+S_L$ denote the coset of a polynomial $p(k)$.
The essence of the polynomial reduction in \cite{WZ2022b} is to find a simple representative for each $[p(k)]_L$.
To this aim, we need to characterize the degree of $L^{\ast}(x(k))$.

Given a nonzero operator $L=\sum\limits_{i=0}^{J}a_i(k)\sigma^i\in \bK[k][\sigma]$ with $a_J(k)\neq 0$,
let
\begin{equation}\label{eq:d and bk}
b_{\ell}(k)=\sum_{j={\ell}}^{J}\binom{j}{\ell}a_{J-j}(k+j-J)\quad
\text{and}
\quad
d=\max_{0\leq \ell \leq J} \{\deg b_\ell(k)-\ell\}.
\end{equation}
For simplicity, we call $d$ in Eq.~\eqref{eq:d and bk} \emph{the degree} of $L$, written as $\deg(L)=d$.
Note that
\begin{equation*}
f(s)=\sum_{\ell=0}^{J}[k^{d+\ell}](b_\ell(k))s^{\underline{\ell}}
\end{equation*}
is a nonzero polynomial in $s$.
Here $[k^{d+\ell}](b_\ell(k))$ denotes the coefficient of $k^{d+\ell}$ in $b_\ell(k)$ and $s^{\underline{\ell}}$ denotes the falling factorial defined by $s^{\underline{\ell}}=s(s-1)\cdots (s-\ell+1)$.
Let
\begin{equation}\label{eq:nonnegative roots}
R_{L}=\{s\in\bN \mid f(s)=0\}.
\end{equation}
Then $L$ is called \emph{nondegenerated} if $R_{L}=\emptyset$, and \emph{degenerated} otherwise.

With the notation above, the degree of $L^{\ast}(x(k))$ is then described as follows.
\begin{lem}\label{lem:degree}\cite[Lemma 2.5]{WZ2022b}
Let $L=\sum\limits_{i=0}^{J}a_i(k)\sigma^i\in\bK[k][\sigma]$ with $a_J(k)\neq 0$ and $d=\deg(L)$ as given by Eq.~\eqref{eq:d and bk}.
Then for any nonzero polynomial $x(k)\in\bK[k]$, we have
\begin{equation*}
\deg L^{\ast}(x(k))\left\{
     \begin{array}{ll}
       < d+\deg x(k), & \hbox{if $L$ is degenerated and $\deg x(k)\in R_{L}$,} \\
       =d+\deg x(k), & \hbox{otherwise.}
     \end{array}
   \right.
\end{equation*}
\end{lem}

By Lemma \ref{lem:degree}, for each $s\in \bN\setminus R_L$, there exists a $q_s(k)\in S_L$ with $\deg q_s(k)=d+s$.
Then using the Euclidean division algorithm, one can write any polynomial $Q(k)\in\bK[k]$ of degree $m$ as
\begin{equation}\label{eq:reduction step1}
 Q(k)=
 \sum_{\substack{0\leq s\leq m-d\\s\notin R_{L}}}u_{s}q_{s}(k)
 +\sum_{\substack{0\leq s\leq m-d\\s\in R_{L}}}u_{s}k^{d+s}+\tilde{q}(k),
\end{equation}
where $u_s\in\bK$ and $\tilde{q}(k)$ is a polynomial with $\deg \tilde{q}(k)<d$.
In particular, when $L$ is nondegenerated, namely, $R_{L}=\emptyset$, we have
\begin{equation}\label{eq:reduction step}
 Q(k)=\sum_{s=0}^{m-d}u_{s}q_{s}(k)+\tilde{q}(k).
\end{equation}
Equation~\eqref{eq:reduction step1} (or Eq.~\eqref{eq:reduction step}) is called a \emph{polynomial reduction} on $Q(k)$ with respect to $L$.
The following theorem follows directly from Eq.~\eqref{eq:reduction step1}.

\begin{theo}
Let $L=\sum_{i=0}^{J}a_i(k)\sigma^i\in \bK[k][\sigma]$ with $a_J(k)\neq 0$, $d=\deg(L)$ and $R_L$ defined by Eq.~\eqref{eq:nonnegative roots}. Then
\begin{equation*}
\bK[k]/S_L=\la[k^{i}]_L\mid i\in\{0,1,2,\ldots,d-1\}\cup (d+R_L)\ra,
\end{equation*}
where $d+R_L=\{d+s\mid s\in R_L\}$.
\end{theo}

In general, one can characterize only the degree but not the structure of $\tilde{q}(k)$.
Nevertheless, when the coefficients of $L$ satisfy additional symmetry condition and $Q(k)=(k-\gamma)^m$ for certain $\gamma\in\bK,\ m\in\bN$, the corresponding $\tilde{q}(k)$ is a linear combination of $(k-\gamma)^j$ with $j$ having the same parity as $m$.

\begin{theo}\label{th:power reduce v2}
Let $L=\sum_{i=0}^{J}a_i(k)\sigma^i\in \bK[k][\sigma]$ with $a_J(k)\neq 0$ and $d=\deg(L)$.
Suppose $L$ is nondegenerated and there exists a $\gamma\in\bK$ such that
\begin{equation}\label{eq:power-reduce-condition}
a_i(\gamma+k)=(-1)^d a_{J-i}(\gamma-k-J),\quad i=0,1,\ldots,\left\lfloor\dfrac{J}{2}\right\rfloor.
\end{equation}
Then for any positive integer $m$, we have
\begin{equation*}
  [(k-\gamma)^m]_L\in\la [(k-\gamma)^i]_L\mid i\equiv m \pmod 2, 0\leq i<d\ra.
\end{equation*}
\end{theo}
When conditions in Theorem \ref{th:power reduce v2} are satisfied, we say $L$ is \emph{power-partible} with respect to $\gamma$.
If $L\in\ann F(k)$ for some holonomic sequence $F(k)$, one may also say $F(k)$ is \emph{power-partible} with respect to $\gamma$.
To prove the theorem, we first recall the following useful observation.
One may refer to \cite[Lemma 3.1]{HouMuZeil2021} for a simple proof.

\begin{lem}\label{lem:HouMuZeil2021}
Let $p(k)\in\bK[k]$ and $\gamma\in\bK$. The following two statements are equivalent:
\begin{enumerate}[(1)]
  \item $p(\gamma+k)=p(\gamma-k)$ ($p(\gamma+k)=-p(\gamma-k)$, respectively);
  \item $p(k)$ is a linear combination of $(k-\gamma)^{2i}$ ($(k-\gamma)^{2i+1}$, respectively), $i\in\bN$.
\end{enumerate}
\end{lem}
\begin{proof}[Proof of Theorem \ref{th:power reduce v2}]
For $s\in\bN$, take
\begin{equation}\label{eq:x}
x_{s}(k)=\left(k-\gamma+\dfrac{J}{2}\right)^s\text{ and }
p_s(k)=L^\ast (x_s(k))=\sum_{i=0}^{J}a_i(k-i)x_s(k-i).
\end{equation}
Then
it is easy to check that
\begin{equation*}x_s(\gamma+k)=(-1)^sx_s(\gamma-k-J).\end{equation*}
Thus
\begin{align*}
p_s(\gamma+k)&=\sum_{i=0}^{J}a_i(\gamma+k-i)x_s(\gamma+k-i)\\
&=\sum_{i=0}^{J}(-1)^d a_{J-i}(\gamma-k+i-J)(-1)^s x_s(\gamma-k+i-J)\\
&=(-1)^{d+s}\sum_{i=0}^{J} a_{J-i}(\gamma-k-(J-i)) x_s(\gamma-k-(J-i))\\
&=(-1)^{d+s}\sum_{i=0}^{J} a_{i}(\gamma-k-i) x_s(\gamma-k-i)\\
&=(-1)^{d+s}p_s(\gamma-k).
\end{align*}
By Lemma \ref{lem:HouMuZeil2021},
$p_s(k)$ is a linear combination of $(k-\gamma)^{2i+1}$ (\emph{resp.} $(k-\gamma)^{2i}$) when $d+s$ is odd (\emph{resp.} even), $i\in\bN$.
Since $L$ is nondegenerated, we know $\deg p_s=d+s$.

Next, we will proceed according to the parity of $d$.
If $d$ is even, then
\begin{align*}
p_{2s}(k)=\sum_{i=0}^{s+d/2}c_{2s,i}(k-\gamma)^{2i}\quad\text{and}\quad
p_{2s+1}(k)=\sum_{i=0}^{s+d/2}c_{2s+1,i}(k-\gamma)^{2i+1},
\end{align*}
with constants $c_{j,i}\in\bK$ and $c_{2s,s+d/2},c_{2s+1,s+d/2}\neq 0$.
When $m$ is even (\emph{resp.} odd), the polynomial reduction of $(k-\gamma)^m$ using $p_{2s}(k)$ (\emph{resp. $p_{2s+1}(k)$}) clearly leads to the conclusion.

If $d$ is odd, then
\begin{align*}
p_{2s}(k)=\sum_{i=0}^{s+(d-1)/2}\tilde{c}_{2s,i}(k-\gamma)^{2i+1} \quad\text{and}\quad
p_{2s+1}(k)=\sum_{i=0}^{s+(d+1)/2}\tilde{c}_{2s+1,i}(k-\gamma)^{2i},
\end{align*}
with constants $\tilde{c}_{j,i}\in\bK$ and $\tilde{c}_{2s,s+(d-1)/2},\tilde{c}_{2s+1,s+(d+1)/2}\neq 0$.
When $m$ is even (\emph{resp.} odd), the polynomial reduction of $(k-\gamma)^m$ using $p_{2s+1}(k)$ (\emph{resp. $p_{2s}(k)$}) also leads to the conclusion.
\end{proof}

From the proof of Theorem \ref{th:power reduce v2}, it is apparent that if we multiply $x_{s}(k)$ in Eq.~\eqref{eq:x} with a nonzero constant $\alpha_s\in\bK$, the rest argument still follows.
With this, we rephrase the above theorem in a more practical form.
\begin{theo}\label{th:power reduce}
Suppose that $L=\sum_{i=0}^{J}a_i(k)\sigma^i$ is power-partible with respect to $\gamma$ and $d=\deg(L)$.
Let $\alpha_s\in\bK^\ast$ and
\begin{equation*}x_{s}(k)=\alpha_s\cdot\left(k-\gamma+\dfrac{J}{2}\right)^s .\end{equation*}
Then for any positive integer $m$, there exist some $u_i,v_j\in\bK$ such that
\begin{equation}\label{eq:power-partible reduction}
(k-\gamma)^m=\sum_{\substack{0\le i<d\\ i\equiv m \pmod 2}}u_i(k-\gamma)^i
+\sum_{\substack{0\le j\le m-d\\ d+j\equiv m \pmod 2}}v_jL^\ast(x_j(k)).
\end{equation}
\end{theo}

\re
We now illustrate that Theorem \ref{th:power reduce v2} is a generalization of both Theorem 3.2 and Theorem 4.2 in \cite{HouMuZeil2021}.
Suppose $t_k$ is hypergeometric and $\dfrac{t_{k+1}}{t_k}=\dfrac{a(k)}{b(k)}$ with $a(k),b(k)\in\bK[k]$.
Then $L=a_0(k)+a_1(k)\sigma\in\ann t_k$ with $a_0(k)=a(k)$ and $a_1(k)=-b(k)$.
Suppose that
\begin{equation*}a(k)=-b(k+\alpha) \quad \text{and}\quad b(\beta+k)=-b(\beta-k)\end{equation*} for some $\alpha,\beta\in\bK$, which is one of the four cases considered in \cite{HouMuZeil2021}.
By Lemma \ref{lem:HouMuZeil2021}, one can see that
$d=\deg(L)=\deg b(k)$ is odd.
Let $\gamma=\beta-\dfrac{\alpha-1}{2}$.
Then it is straightforward to check that
\begin{equation*}
a_0(\gamma+k)
=-b\left(\beta+\dfrac{\alpha+1}{2}+k\right)
=b\left(\beta-\dfrac{\alpha+1}{2}-k\right)
=(-1)^d a_1(\gamma-k-J),
\end{equation*}
where $J=1$ is the order of $L$.
That is, Eq.~\eqref{eq:power-reduce-condition} holds.
By similar discussions, it can be easily checked that Eq.~\eqref{eq:power-reduce-condition} always holds whenever $a(k)=\pm b(k+\alpha)$ and $b(\beta+k)=\pm b(\beta-k)$.

\section{Congruences for Ap\'ery numbers}

In this section, we show that the Ap\'ery numbers $A_k$ are power-partible, and then use this property to derive new series of congruences for them.

Recall that the Ap\'ery numbers are defined as
\begin{equation*}
A_n=\sum_{k=0}^{n}\binom{n}{k}^2\binom{n+k}{k}^2.
\end{equation*}
By Zeilberger's algorithm, we find that
\begin{equation}\label{eq:L Apery1}
L=a_2(k)\sigma^2+a_1(k)\sigma+a_0(k)\in \ann A_k,
\end{equation}
where
\begin{equation*}a_2(k)=(k+2)^3,a_1(k)=-(2k+3)(17k^2+51k+39),a_0(k)=(k+1)^3.\end{equation*}
It is easy to check that $d=\deg(L)=3$, $L$ is nondegenerated and
\begin{equation*}a_0(\gamma+k)=(-1)^da_2(\gamma-k-2),a_1(\gamma+k)=(-1)^da_1(\gamma-k-2),\end{equation*}
where $\gamma=-\dfrac{1}{2}$.
That is, the Ap\'ery numbers $A_k$ are power-partible with respect to $\gamma=-\dfrac{1}{2}$.
By Theorem \ref{th:power reduce v2} we know
\begin{equation}\label{eq:reT(2k+1)}
[(2k+1)^{2r+1}]_L\in \la[(2k+1)]_L\ra.
\end{equation}

The following lemma follows from Theorem 1.1 (ii) in \cite{Sun2012} and the fact that the denominator of the Bernoulli number $B_{2k}$ is square-free \cite[p.116]{Hardy--Wright}.
\begin{lem}[Sun \cite{Sun2012}]\label{lem:Sun-Apery}
Let $p>3$ be a prime. Then
\begin{equation*}
\sum_{k=0}^{p-1}(2k+1)A_k\equiv p \pmod {p^3}.
\end{equation*}
\end{lem}
Lemma \ref{lem:Sun-Apery} and Eq.~\eqref{eq:reT(2k+1)} motivate the discovery of the following congruences.
\begin{theo}\label{th:Ak(2k+1)^{2r+1}}
For each $r\in\bN$,
there exists an integer $\tilde{c}_r$ such that
\begin{equation}\label{1q:(2k+1)^{2r+1}}
\sum_{k=0}^{p-1}(2k+1)^{2r+1}A_k\equiv \tilde{c}_r p \pmod {p^3}
\end{equation}
holds for any prime $p>3$.
\end{theo}
To prove Theorem \ref{th:Ak(2k+1)^{2r+1}}, we first discuss the arithmetic properties of $L^{\ast}(x(k)) A_k$.
\begin{lem}\label{th:Apery cong}
Let $L$ be given as in Eq.~\eqref{eq:L Apery1} and $n$ a positive integer. Then for any polynomial $x(k)\in\bK[k]$, we have
\begin{equation}\label{eq: Apery D1}
\sum_{k=0}^{n-1} L^{\ast}(x(k)) A_k
=n^3\left(x(n-1)A_{n-1}-x(n-2)A_n\right).
\end{equation}
Here $L^{\ast}$ is the adjoint of $L$.
\end{lem}
\begin{proof}
By Eq.~\eqref{eq:rec ann} and the fact $u_0(0)A_0+u_1(0)A_1=0$, we have
\begin{equation}\label{eq:Apery Delta1}
\sum_{k=0}^{n-1}L^{\ast}(x(k)) A_k
=
-\left(u_0(n)A_n+u_1(n)A_{n+1}\right),
\end{equation}
where
$u_0(n)=n^3x(n-2)-(2 n+1) (17 n^2+17n+5) x(n-1)$ and
$u_1(n)=(n+1)^3 x(n-1)$.
As $L\in\ann A_k$, it is clear that for any $n\geq 1$,
\begin{equation}\label{eq:shift Apery1}
(n+1)^3A_{n+1}=(2n+1)(17n^2+17n+5)A_n-n^3A_{n-1}.
\end{equation}
Substituting Eq.~\eqref{eq:shift Apery1} into Eq.~\eqref{eq:Apery Delta1}, we derive Eq.~\eqref{eq: Apery D1}.
\end{proof}

\begin{coro}\label{coro:Apery cong1}
Let $L$ be given as in Eq.~\eqref{eq:L Apery1} and $n$ a positive integer. Then for any polynomial $x(k)\in\bZ[k]$, we have
\begin{equation*}\label{eq:Apery cong1}
\sum_{k=0}^{n-1} L^{\ast}(x(k))A_k \equiv 0 \pmod {n^3}.
\end{equation*}
\end{coro}
Recall that $L$ in Eq.~\eqref{eq:L Apery1} is power-partible with respect to $\gamma=-\dfrac{1}{2}$.
Let
\begin{equation}\label{eq:xs1}
x_s(k)=2^{s+1}\left(k-\gamma+\dfrac{2}{2}\right)^s=2(2k+3)^s,\quad s\in\bN.
\end{equation}
Then by the proof of Theorem \ref{th:power reduce v2}, $L^\ast(x_s(k))$ is a linear combination of $(2k+1)^{i}$ with $i\equiv s+3{\pmod 2}$.
We show that the coefficients in the combination are all integers, and that the coefficients are all divisible by $8$ when $s$ is even.
\begin{lem}\label{lm:Structure}
Let $L$ be given by Eq.~\eqref{eq:L Apery1} and $x_s(k)$ given by Eq.~\eqref{eq:xs1}.
Then
\begin{equation}\label{eq:Explicit form of L*(x_s(k))1}
L^\ast(x_s(k))=-8(2k+1)^{s+3}+\sum_{j=1}^{\left\lfloor \frac{s+3}{2}\right\rfloor}e_j(2k+1)^{s+3-2j},
\end{equation}
where $e_j\in\bZ$ for all $j=1,2,\ldots,\left\lfloor \dfrac{s+3}{2}\right\rfloor$.
If $s$ is even, we further have
\begin{equation}\label{eq:Explicit form of L*(x_s(k))2}
8\mid e_j,\quad \forall j=1,2,\ldots,s/2+1.
\end{equation}
\end{lem}
\begin{proof}
For simplicity, let $\ell=2k+1$. By the definition of $L^\ast$, we have
\begin{align*}
 &L^\ast(x_s(k))=\sum_{i=0}^{2}a_{i}(k-i)x_s(k-i)\\
=&2(k+1)^3(2k+3)^s-2(17 k^2+17 k+5)(2k+1)^{s+1}+2k^3(2k-1)^s\\
=&\dfrac{(\ell+1)^3}{4}(\ell+2)^s-\dfrac{1}{2}(17\ell^2+3)\ell^{s+1}
+\dfrac{(\ell-1)^3}{4}(\ell-2)^s\\
=&\dfrac{(\ell+1)^3}{4}\sum_{j=0}^{s}\binom{s}{j}2^{j}\ell^{s-j}
+\dfrac{(\ell-1)^3}{4}\sum_{j=0}^{s}\binom{s}{j}(-2)^{j}\ell^{s-j}
-\dfrac{\ell^s(17\ell^3+3\ell)}{2}\\
=&\dfrac{(\ell^3+3\ell)}{2}\sum_{\substack{j=0\\ j\text{ even}}}^{s}\binom{s}{j}2^{j}\ell^{s-j}
+\dfrac{(3\ell^2+1)}{2}\sum_{\substack{j=0\\ j\text{ odd}}}^{s}\binom{s}{j}2^{j}\ell^{s-j}
-\dfrac{\ell^s(17\ell^3+3\ell)}{2}
\end{align*}
\begin{align}\label{eq:r1}
=&-8\ell^{s+3}+\sum_{\substack{j=2\\ j\text{ even}}}^{s}\binom{s}{j}2^{j-1}\ell^{s-j}(\ell^3+3\ell)+\sum_{\substack{j=1\\ j\text{ odd}}}^{s}\binom{s}{j}2^{j-1}\ell^{s-j}(3\ell^2+1)\\
=&-8\ell^{s+3}+\sum_{j=1}^{\left\lfloor \frac{s+3}{2}\right\rfloor}e_j\ell^{s+3-2j}\nonumber
\end{align}
where $e_j\in\bZ$ for all $j=1,2,\ldots,\left\lfloor \dfrac{s+3}{2}\right\rfloor$.
At this stage, Eq.~\eqref{eq:Explicit form of L*(x_s(k))1} has been proved.
Now we consider the case when $s$ is even.
It is direct to check that
\[L^\ast(x_{0}(k))=-8\ell^3\text{ and }
L^\ast(x_{2}(k))=-8 (\ell^5-\ell^3-\ell).
\]
Next, we assume $s$ is even and $s>2$.
By Eq.~\eqref{eq:r1}, we only need to consider the terms in the summations for $j<4$ since $8\mid 2^{j-1}$ for $j\ge 4$.
Since
\begin{align*}
& \binom{s}{2}2\ell^{s-2}(\ell^3+3\ell)+
  \binom{s}{1} \ell^{s-1}(3\ell^2+1)+
  \binom{s}{3}4\ell^{s-3}(3\ell^2+1)\\
=& \left(2\binom{s}{2}+3\binom{s}{1}\right)\ell^{s+1}
+\left(6\binom{s}{2}+\binom{s}{1}+12\binom{s}{3}\right)
\ell^{s-1}
+4\binom{s}{3}\ell^{s-3}\\
=&s(s+2)\ell^{s+1}+(2s^2(s-1)-s(s-2))\ell^{s-1}+\dfrac{2s(s-1)(s-2)}{3}\ell^{s-3}.
\end{align*}
It is then straightforward to check that the above three coefficients are all multiples of eight.
This completes the proof of Eq.~\eqref{eq:Explicit form of L*(x_s(k))2}.
\end{proof}

\begin{proof}[Proof of Theorem \ref{th:Ak(2k+1)^{2r+1}}]
Let $L$ be given by Eq.~\eqref{eq:L Apery1} and $x_s(k)=2(2k+3)^s$ as in Eq.~\eqref{eq:xs1}.
By Eq.~\eqref{eq:Explicit form of L*(x_s(k))2} and the expression for $L^\ast(x_{s}(k))$ in Eq.~\eqref{eq:Explicit form of L*(x_s(k))1}, the power-partible reduction on $(2k+1)^{2r+1}$ reveals that
\begin{equation}\label{q:(2k+1)^{2+1}decomp}
(2k+1)^{2r+1}=\sum_{s=0}^{r-1}
v_s\left(\dfrac{1}{8}L^\ast(x_{2s}(k))\right)+\tilde{c}_r(2k+1),
\end{equation}
for some $\tilde{c}_r,v_s\in\bZ$.
Multiplying both sides of Eq.~\eqref{q:(2k+1)^{2+1}decomp} with $A_k$ and then
summing over $k$ from $0$ to $p-1$, we obtain from Corollary \ref{coro:Apery cong1} that
\begin{equation*}\sum_{k=0}^{p-1}(2k+1)^{2r+1}A_k\equiv \tilde{c}_r\sum_{k=0}^{p-1}(2k+1)A_k \pmod {p^3}\end{equation*}
since $p>3$ is a prime. Equation \eqref{1q:(2k+1)^{2r+1}} then follows from Lemma \ref{lem:Sun-Apery}.
\end{proof}

The proof also presents an algorithmic process to determine the $\tilde{c}_r$ for each $r\in\bN^\ast$.
For example, when $r=1$, since $(2 k + 1)^3=-\dfrac{1}{8}L^{\ast}(x_0(k))$, we obtain $\tilde{c}_1=0$ and
\begin{equation*}
\sum_{k=0}^{p-1}(2k+1)^3A_k\equiv 0  \pmod {p^3}.
\end{equation*}
When $r=2$, we have
$(2 k + 1)^5=(2k+1)-\dfrac{1}{8}L^{\ast}(x_0(k))-\dfrac{1}{8}L^{\ast}(x_2(k))$.
Thus $\tilde{c}_2=1$ and
\begin{equation*}
\sum_{k=0}^{p-1}(2k+1)^5 A_k\equiv p \pmod {p^3}.
\end{equation*}

The above argument also applies to the alternating case, namely,
Eq.~\eqref{eq: Sun-Conj-Alternating} can also be proved using the reduction method.
Let $F(k)=(-1)^kA_k$, Zeilberger's algorithm leads to
\begin{equation*}\label{eq:L Apery}
\tilde{L}=\tilde{a}_2(k)\sigma^2+\tilde{a}_1(k)\sigma+\tilde{a}_0(k)\in \ann F(k),
\end{equation*}
where
\begin{equation*}\tilde{a}_2(k)=(k+2)^3,\tilde{a}_1(k)=(2k+3)(17k^2+51k+39),\tilde{a}_0(k)=(k+1)^3.\end{equation*}
It can be easily checked that $\deg(\tilde{L})=3$, $\tilde{L}$ is nondegenerated and power-partible with respect to $\gamma=-\dfrac{1}{2}$.
By Theorem \ref{th:power reduce v2} we know
\begin{equation*}\label{eq:-re(2k+1)}
[(2k+1)^{2r+1}]_{\tilde{L}}\in \la[(2k+1)]_{\tilde{L}}\ra.
\end{equation*}
By similar discussions to Lemmas \ref{th:Apery cong} and \ref{lm:Structure}, we can obtain
\begin{equation}\label{eq: Apery -1}
\sum_{k=0}^{n-1} \tilde{L}^{\ast}(x(k))(-1)^kA_k
=n^3\left(x(n-1)F(n-1)-x(n-2)F(n)\right),
\end{equation}
for any polynomial $x(k)\in\bK[k]$
and
\begin{equation}\label{eq:-Explicit form of L*(x_s(k))}
\tilde{L}^\ast(x_s(k))=9(2k+1)^{s+3}+\sum_{j=1}^{\left\lfloor \frac{s+3}{2}\right\rfloor}f_j(2k+1)^{s+3-2j},
\end{equation}
where $f_j\in\bZ$ for all $j=1,2,\ldots,\left\lfloor\dfrac{s+3}{2}\right\rfloor$.
By the reduction on $(2k+1)^{2r+1}$ with respect to $\tilde{L}^\ast(x_s(k))$, one can then prove Eq.~\eqref{eq: Sun-Conj-Alternating} by using Eqs.~\eqref{eq:(2k+1)Sun-Guo-Zeng}, \eqref{eq: Apery -1} and \eqref{eq:-Explicit form of L*(x_s(k))}.

\section*{Acknowledgments}
This work was supported by the National Natural Science Foundation of China (No. 12101449, 12271511, 12271403) and the Natural Science Foundation of Tianjin, China (No. 22JCQNJC00440).
We also would like to express our sincere gratitude to the anonymous reviewer for his valuable comments, which have greatly improved this paper.

\end{document}